\documentclass[12pt]{article}
\usepackage{amssymb}
\usepackage{mathrsfs}
\usepackage{pifont}
\usepackage{tipa}
\usepackage{amsmath}
\usepackage{amssymb,amsmath}
\usepackage{amsfonts}
\catcode`@=11
 \long\def\@makefntext#1{\noindent #1}

\begin{document}
\title{\large{\textbf{Classification of simple weight modules for the
 Neveu-Schwarz algebra with a finite-dimensional weight space}} }
\author{Xiufu Zhang $^{1}$, Zhangsheng Xia $^2$\\
{\scriptsize 1. School of Mathematical Sciences, Xuzhou Normal
University,
Xuzhou 221116,  China}\\
{\scriptsize 2. School of Sciences, Hubei University for
Nationalities, Enshi 445000, China}}
\date{}
\maketitle \footnotetext{\footnotesize* Supported by the National
Natural Science Foundation of China (No. 10931006).}
\footnotetext{\footnotesize** Email: xfzhang@xznu.edu.cn;\quad
xzsw8577@163.com}

\numberwithin{equation}{section}

\begin{abstract}
We show that the support of a simple weight module over the
Neveu-Schwarz algebra, which has an infinite-dimensional weight
space, coincides with the weight lattice and that all non-trivial
weight spaces of such module are infinite-dimensional. As a
corollary we obtain that every simple weight module over the
Neveu-Schwarz algebra, having a non-trivial finite-dimensional
weight space, is a Harish-Chandra module (and hence is either a
highest or lowest weight module, or else a module of the
intermediate series). This result generalizes a theorem which was
originally given on the Virasoro algebra. \vspace{2mm}
\vspace{2mm}
\end{abstract}

\vskip 3mm\noindent$1.\ \mathbf{Introduction}$

\vskip 3mm It is well known that the Virasoro algebra $Vir$ plays a
fundamental role in two-dimensional conformal quantum field theory.
As the super-generalization, there are two super-Virasoro algebras
called the \emph{Neveu-Schwarz algebra} and the \emph{Ramond
algebra} corresponding to $N=\frac{1}{2}$ and $N=1$ super-conformal
field theory respectively. Let $\theta=\frac{1}{2}$ or $0$ which
corresponds to the Neveu-Schwarz case or the Ramond case
respectively. The super-Virasoro algebra $SVir(\theta)$ is the Lie
superalgebra $SVir(\theta)=SVir_{\bar{0}}\oplus SVir_{\bar{1,}}$
where $SVir_{\bar{0}}$ has a basis $\{L_n, c\mid n\in \mathbb{Z}\}$
and $SVir_{\bar{1}}$ has a basis $\{G_r\mid r\in
\theta+\mathbb{Z}\},$ with the commutation relations
\begin{align*}
&[L_m,L_n]=(n-m)L_{n+m}+\delta_{m+n,0}\frac{m^{3}-m}{12}c,\\
&[L_m,G_r]=(r-\frac{m}{2})G_{m+r},\\
&[G_r,G_s]=2L_{r+s}+\frac{4r^2-1}{12}\delta_{r+s,0}c,\\
&[SVir_{\bar{0}},c]=0=[SVir_{\bar{1}},c],
\end{align*}
for $m, n\in\mathbb{Z}, r, s \in \theta+\mathbb{Z}.$

The subalgebra $\mathfrak{h}=\mathbb{C}L_0\oplus\mathbb{C}c$ is
called the \emph{Cartan subalgebra} of $SVir(\theta).$ An
$\mathfrak{h}$-diagonalizable  $SVir(\theta)$-module is usually
called a \emph{weight module}. If $M$ is a weight module, then $M$
can be written as a direct sum of its weight spaces, $M=\bigoplus
M_{\lambda},$ where $M_{\lambda}=\{v\in M|L_0v=\lambda(L_0) v,
cv=\lambda(c) v\}.$ We call $\{\lambda|M_{\lambda}\neq0\}$ the
$support$ of $M$ and is denoted by $supp(M).$ A weight module is
called \emph{Harish-Chandra module} if each weight space is
finite-dimensional.

In [1], [2], [3], [7],  the unitary modules, the Verma modules, the
Fock modules and the Harish-Chadra modules over $SVir(\theta)$ are
studied. In this paper, we try to give the classification of simple
weight modules for the Neveu-Schwarz algebra with a
finite-dimensional weight space. Our main result generalizes a
theorem which was originally given on the Virasoro algebra in [6].

For convenience, we denote the Neveu-Schwarz algebra by $NS$ instead
of $SVir(\frac{1}{2})$ and we define
$$E_{k}=\left\{\begin{array}
{l@{\quad\quad}l} L_k, & \mathrm{if}\ \ k\in \mathbb{Z}, \\
G_{k}, &\mathrm{if}\ \ k\in \frac{1}{2}+\mathbb{Z}.
\end{array}\right.$$

Suppose $M$ is a simple weight $NS-$module, then $c$ acts on $M$ by
a scalar and $M$ can be written as a direct sum of its weight
spaces, $M=\bigoplus M_{\lambda},$ where $M_{\lambda}=\{v\in
M|L_0v=\lambda v\}.$ Obviously, if $\lambda\in \mathrm{supp}(M),$
then $\mathrm{supp}(M)\subseteq
\{\lambda+k|k\in\frac{\mathbb{Z}}{2}\},$ the weight lattice.  Two
elements $i,j\in\{\lambda+k|k\in\frac{\mathbb{Z}}{2}\}$  are called
$adjacent$ if $|i-j|=\frac{1}{2},$ otherwise, $unadjacent.$

\vskip 3mm In this paper, our main result is the following theorem:

\vskip 3mm\noindent$\mathbf{Theorem\ 1.}$ Let $M$ be a simple weight
$NS-$module. Assume that there exists $\lambda\in \mathbb{C}$ such
that $\mathrm{dim}M_{\lambda}=\infty.$ Then
$supp(M)=\lambda+\frac{\mathbb{Z}}{2}$ and for each $k\in
\frac{\mathbb{Z}}{2}$ we have $\mathrm{dim}(M_{\lambda+k})=\infty.$

\vskip 3mm In [7], Y. Su proved the following result which
generalizes a theorem originally given as a conjecture by Kac in [4]
on the Virasoro algebra and proved by Mathieu in [5].

\vskip 3mm\noindent$\mathbf{Theorem\ 2.}$ A Harish-Chandra module
over $SVir(\theta)$ is either a highest or lowest weight module, or
else a module of the intermediate series.

\vskip 3mm By Theorem 1 and 2, we get the following corollaries
immediately:

\vskip 3mm\noindent$\mathbf{Corrollary\ 3.}$ Let $M$ be a simple
weight $NS-$module. Assume that there exists $\lambda\in \mathbb{C}$
such that $0<\mathrm{dim}M_{\lambda}<\infty.$ Then $M$ is a
Harish-Chandra module. Consequently, $M$ is either a highest or
lowest weight module, or else a module of the intermediate series.

\vskip 3mm A weight $NS$-module, $M,$ is called \emph{mixed} module
if there exist $\lambda\in \mathbb{C}$ and
$k\in\frac{\mathbb{Z}}{2}$ such that
$\mathrm{dim}M_{\lambda}=\infty$ and
$\mathrm{dim}M_{\lambda+k}<\infty.$

\vskip 3mm\noindent$\mathbf{Corrollary\ 4.}$ There are no simple
mixed $NS$-modules.

\vskip 3mm\noindent$2.\ \mathbf{Proof\ of \ the\ Theorem\ 1}$

\vskip 3mm Noting that
$\{G_{-\frac{3}{2}},G_{-\frac{1}{2}},G_{\frac{1}{2}},G_{\frac{3}{2}}\}$
is a set of generators of $NS$, we have the following fact
immediately:

\vskip 3mm\noindent$\mathbf{Principal\ \ Fact}$: Assume that there
exists $\mu\in \mathbb{C}$ and  $0\neq v\in M_\mu,$ such that
$G_{\frac{1}{2}}v=G_{\frac{3}{2}}v=0$ or
$G_{-\frac{1}{2}}v=G_{-\frac{3}{2}}v=0$. Then $M$ is a
Harish-Chandra $NS$-module.

\vskip 3mm\noindent{\bf{Lemma 1.}} If there are
$\mu,\mu^{'}\in\{\lambda+k|k\in\frac{\mathbb{Z}}{2}\}$ such that
$\mathrm{dim}M_\mu<\infty$ and $\mathrm{dim}M_{\mu^{'}}<\infty,$
then $\mu,\mu^{'}$ are adjacent.

\vskip 3mm\noindent{\bf{Proof.}} Suppose that there exist two
unadjacent elements in $\{\lambda+k|k\in\frac{\mathbb{Z}}{2}\}$
correspond to finite-dimensional weight spaces or trivial vector
spaces in $M.$ Without loss of generality, we may assume that
$\mathrm{dim}M_{\lambda+\frac{1}{2}}<\infty$ and
$\mathrm{dim}M_{\lambda+k}<\infty$, where
$k\in\{\frac{3}{2},2,\frac{5}{2},3,\cdots\}$. Let $V$ denote the
intersection of the kernels of the linear maps:
$$
G_{\frac{1}{2}}: M_{\lambda}\rightarrow M_{\lambda+\frac{1}{2}},$$
and $$E_{k}: M_{\lambda}\rightarrow M_{\lambda+k}.
$$
Since $\mathrm{dim}M_{\lambda}=\infty,
\mathrm{dim}M_{\lambda+\frac{1}{2}}<\infty$, we know that the kernel
$kerG_{\frac{1}{2}}$ of $G_{\frac{1}{2}}:M_\lambda\rightarrow
M_{\lambda+\frac{1}{2}}$ is infinite dimensional. Since
$\mathrm{dim}M_{\lambda+k}<\infty$, we also have that the kernel of
$E_{k}|_{kerG_{\frac{1}{2}}}:kerG_{\frac{1}{2}}\rightarrow
M_{\lambda+k}$ is infinite dimensional. That is
$\mathrm{dim}V=\infty.$

Since
$$
[G_{\frac{1}{2}},E_l]=\left\{\begin{array}
{l@{\quad\quad}l} 2L_{\frac{1}{2}+l}, & \mathrm{if}\ \ l\in \frac{1}{2}+\mathbb{N}, \\
(\frac{l}{2}-\frac{1}{2})G_{\frac{1}{2}+l}, &\mathrm{if}\ \ l\in
\mathbb{Z}_{+}\setminus\{1\}
\end{array}\right.
$$
is nonzero, we get that $$E_lV=0,\ \forall l=\frac{1}{2}, k,
k+\frac{1}{2}, k+1, \cdots  \eqno(1)$$ If $k=\frac{3}{2}$ then there
exists $0\neq v\in V$ such that $G_{\frac{3}{2}}v=0,$ and $M$ would
be a Harish-Chandra module by the Principal Fact, a contradiction.
Thus $k>\frac{3}{2}$ and $\mathrm{dim}G_{\frac{3}{2}}V=\infty.$
Since $\mathrm{dim}M_{\lambda+\frac{1}{2}}<\infty,$ there exists
$0\neq w\in G_{\frac{3}{2}}V$ such that  $L_{-1}w=0.$ Suppose
$w=G_{\frac{3}{2}}u$ for some $u\in V.$ For all $l\geq
k>\frac{3}{2},$ using (1) we have
$$E_lw=E_lG_{\frac{3}{2}}u=\left\{\begin{array}{rl}
-G_{\frac{3}{2}}E_lu+[E_l,G_{\frac{3}{2}}]u=0, & \mathrm{if}\ \ k\in \frac{1}{2}+\mathbb{Z}_{+},\\
G_{\frac{3}{2}}E_lu+[E_l,G_{\frac{3}{2}}]u=0, &\mathrm{if}\ \ k\in
\mathbb{Z}_{+}.
\end{array}\right.$$
Hence $$E_lw=0,\ \forall l=-1,k,k+\frac{1}{2},k+1.\cdots$$ Since
$[L_{-1},E_{l}]\neq0$ for all $l\geq\frac{3}{2}$, inductively, we
get $$E_lw=0,\ \forall l=\frac{1}{2},1,\frac{3}{2},2,\cdots.$$ Hence
$M$ is a Harish-Chandra module by the Principal Fact. A
contradiction. Then the lemma follows.
 \hfill$\Box$

By Lemma 1, we see that there exist at most two elements in
$\{\lambda+k|k\in\frac{\mathbb{Z}}{2}\}$ which
 correspond to
finite-dimensional weight spaces or trivial vector spaces in $M$.
Moreover, if there are two, then they are adjacent.

\vskip 3mm\noindent{\bf{Lemma 2.}} (i) Let $0\neq v\in M$ be such
that $G_{\frac{1}{2}}v=0.$ Then
$$
(\frac{1}{2}L_1G_{\frac{1}{2}}-G_{\frac{3}{2}})G_{\frac{3}{2}}v=0.
$$

(ii) Let $0\neq w\in M$ be such that $G_{-\frac{1}{2}}w=0.$ Then
$$
(\frac{1}{2}L_{-1}G_{-\frac{1}{2}}+G_{-\frac{3}{2}})G_{-\frac{3}{2}}w=0.
$$

\vskip 3mm\noindent{\bf{Proof.}} Note that
$L_1v=G_{\frac{1}{2}}G_{\frac{1}{2}}v=0,
L_{-1}w=G_{-\frac{1}{2}}G_{-\frac{1}{2}}w=0,$ we can easily check by
a direct calculation that Lemma 2 holds. \hfill$\Box$

\vskip 3mm\noindent{\bf{Lemma 3.}} Let $M$ be a simple weight
$NS$-module satisfying  $\mathrm{dim}M_{\mu}<\infty$ and
$\mathrm{dim}M_{\mu+\frac{1}{2}}<\infty.$ Then $\mu\in\{-1,
\frac{1}{2}\}.$

\vskip 3mm\noindent{\bf{Proof.}} Let $V$ be the kernel of
$G_{\frac{1}{2}}: M_{\mu-\frac{1}{2}}\rightarrow M_{\mu}$. Since
$\mathrm{dim}M_{\mu-\frac{1}{2}}=\infty$ and
$\mathrm{dim}M_{\mu}<\infty,$ we see that $\mathrm{dim}V=\infty.$
For any $0\neq v \in V,$ consider the element $G_{\frac{3}{2}}v.$ By
the Principal Fact, $G_{\frac{3}{2}}v=0$ would imply that $M$ is a
Harish-Chandra module, a contradiction. Hence
$G_{\frac{3}{2}}v\neq0,$ in particular,
$\mathrm{dim}G_{\frac{3}{2}}V=\infty.$  This implies that there
exists $w\in G_{\frac{3}{2}}V$ such that $w\neq0$ and
$G_{-\frac{1}{2}}w=0$ since
$\mathrm{dim}M_{\mu+\frac{1}{2}}<\infty.$ From Lemma 2 we have
$(\frac{1}{2}L_1G_{\frac{1}{2}}-G_{\frac{3}{2}})w=0.$ In particular,
we have
$L_{-1}G_{-\frac{1}{2}}(\frac{1}{2}L_1G_{\frac{1}{2}}-G_{\frac{3}{2}})w=0.$
By a direct calculation we obtain
$$
L_{-1}G_{-\frac{1}{2}}(\frac{1}{2}L_1G_{\frac{1}{2}}-G_{\frac{3}{2}})
\equiv 2L_0^{2}-3L_0 \ \mathrm{mod}\ U(NS)G_{-\frac{1}{2}}.
$$
But $w\in M_{\mu+1,}$ which means $L_0w=(\mu+1)w$ and hence
$$
2(\mu+1)^{2}-3(\mu+1)=0.
$$
So $\mu\in \{-1,\frac{1}{2}\}.$ \hfill$\Box$

From Lemma 3 we know that if $M$ has two finite dimension weight
spaces, then they must be $M_{-1},$  $M_{-\frac{1}{2}}$ or
$M_{\frac{1}{2}}, M_{1}.$

\vskip 3mm\noindent{\bf{Lemma 4.}} $\mu\neq-1$ and
$\mu\neq\frac{1}{2}$.

\vskip 3mm\noindent{\bf{Proof.}} By the symmetry, we need only to
prove that $\mu\neq \frac{1}{2}$. Let $V$ be the kernel of
$G_{\frac{1}{2}}: M_{0}\rightarrow M_{\frac{1}{2}}.$ Then
$\mathrm{dim}V=\infty.$ For $v\in V,$ using
$G_{\frac{1}{2}}v=L_0v=0,$ we have
$$G_{\frac{1}{2}}G_{-\frac{1}{2}}v=2L_0v-G_{-\frac{1}{2}}G_{\frac{1}{2}}v=0. \eqno(2)$$
If $G_{-\frac{1}{2}}V$ were infinite-dimensional, there would exists
$0\neq w\in G_{-\frac{1}{2}}V$ such that $G_{\frac{1}{2}}w=0$ (by
(2)) and $G_{\frac{3}{2}}w=0$ (since $\mathrm{dim}M_{1}<\infty$).
Then the Principal Fact would then imply that $M$ is a
Harish-Chandra module, a contradiction. Hence
$$\mathrm{dim}G_{-\frac{1}{2}}V<\infty.$$ This means that the kernel
$W$ of the linear map $$G_{-\frac{1}{2}}: V\rightarrow
M_{-\frac{1}{2}}$$ is infinite-dimensional. For every $x\in W$ we
have
$$G_{\frac{1}{2}}G_{-\frac{3}{2}}x=2L_{-1}x-G_{-\frac{3}{2}}G_{\frac{1}{2}}x=0. \eqno(3)$$
If there exists $0\neq x\in W$ such that $G_{-\frac{3}{2}}x=0$, then
we would have $G_{-\frac{3}{2}}x=G_{-\frac{1}{2}}x=0$ and the
Principal Fact would imply that $M$ is a Harish-Chandra module, a
contradiction. Thus $$\mathrm{dim}G_{-\frac{3}{2}}W=\infty.$$ Let
$H$ denote the kernel of the linear map $L_2:
G_{-\frac{3}{2}}W\rightarrow M_{\frac{1}{2}}.$ Since
$\mathrm{dim}G_{-\frac{3}{2}}W=\infty$ and
$\mathrm{dim}M_{\frac{1}{2}}<\infty,$ we have
$\mathrm{dim}H=\infty.$ For every $y\in H,$ we also have
$G_{\frac{1}{2}}y=0$ by (3), implying by induction that $E_{k}H=0$
for all $k=\frac{1}{2},1,2,\frac{5}{2},3,\frac{7}{2},4,\cdots.$ If
$G_{\frac{3}{2}}h=0$ for some $0\neq h\in H$, then the Principal
Fact implies that $M$ is a Harish-Chandra module, a contradiction.
Hence $$\mathrm{dim}G_{\frac{3}{2}}H=\infty.$$ For every $h\in H$
and $k\geq 2,$ we have
$$
E_kG_{\frac{3}{2}}h=\left\{\begin{array}
{l@{\quad\quad}l} 2L_{\frac{3}{2}+k}h-G_{\frac{3}{2}}G_{k}h=0, & \mathrm{if}\ \ k\in \frac{1}{2}+\mathbb{Z}_{+},\\
(\frac{3}{2}-\frac{k}{2})G_{\frac{3}{2}+k}h+G_{\frac{3}{2}}L_{k}h=0,
&\mathrm{if}\ \ k\in \mathbb{Z}_{+}.
\end{array}\right.
$$
Hence $$E_kG_{\frac{3}{2}}h=0,\  \forall k=2, \frac{5}{2}, 3,
\frac{7}{2}, 4, \frac{9}{2}, \cdots. \eqno(4)$$ Let, finally, $K$
denote the infinite-dimensional kernel of the linear map
$$G_{\frac{1}{2}}: G_{\frac{3}{2}}H\rightarrow M_{\frac{1}{2}}.$$
If $G_{\frac{3}{2}}z=0$ for some $0\neq z\in K,$ then the Principal
Fact implies that $M$ is a Harish-Chandra module, a contradiction.
Hence $G_{\frac{3}{2}}z\neq 0, \forall 0\neq z\in K.$ For every
$z\in K$ and $k\geq 2$, by (4), we have
$\mathrm{dim}G_{\frac{3}{2}}K=\infty$ and
$$
E_kG_{\frac{3}{2}}z=\left\{\begin{array}
{l@{\quad\quad}l} 2L_{\frac{3}{2}+k}z-G_{\frac{3}{2}}G_{k}z=0, & \mathrm{if}\ \ k\in \frac{1}{2}+\mathbb{Z}_{+},\\
(\frac{3}{2}-\frac{k}{2})G_{\frac{3}{2}+k}z+G_{\frac{3}{2}}L_{k}z=0.
&\mathrm{if}\ \ k\in \mathbb{Z}_{+}.
\end{array}\right.
$$
Hence $E_kG_{\frac{3}{2}}K=0$ for all $k\geq2.$ At the same time,
since $\mathrm{dim}G_{\frac{3}{2}}K=\infty$ and
$\mathrm{dim}M_1<\infty,$ we can choose some $0\neq t\in
G_{\frac{3}{2}}K$ such that $G_{-\frac{1}{2}}t=0,$ by induction, we
get $E_it=0$ for all $i>0$ and thus $M$ is a Harish-Chandra module
by the Principal Fact. This last contradiction completes the proof
of lemma 4.\hfill$\Box$

\vskip3mm By lemma 1, Lemma 3 and Lemma 4, we get the following
lemma immediately:

\vskip3mm\noindent{\bf{Lemma 5.}} There is at most one element
$\mu\in \{\lambda+k|k\in\frac{\mathbb{Z}}{2}\}$ such that
$\mathrm{dim}M_{\mu}<\infty.$

\vskip3mm Now the proof of Theorem 1 follows from the following
Lemma:

\vskip3mm\noindent{\bf{Lemma 6.}} There is no $\mu\in
\{\lambda+k|k\in\frac{\mathbb{Z}}{2}\}$ such that
$\mathrm{dim}M_{\mu}<\infty.$

\vskip3mm\noindent{\bf{Proof.}} Suppose that
$\mathrm{dim}M_{\mu}<\infty$ and $\mathrm{dim}M_{\nu}=\infty$ for
all $\mu\neq\nu\in\{\lambda+k|k\in\frac{\mathbb{Z}}{2}\}.$ Define
$$V=ker(G_{\frac{1}{2}}:M_{\mu-\frac{1}{2}}\rightarrow M_{\mu})\cap
ker(G_{\frac{1}{2}}G_{-\frac{3}{2}}G_{\frac{3}{2}}:M_{\mu-\frac{1}{2}}\rightarrow
M_{\mu})\cap
ker(L_{-1}G_{\frac{3}{2}}:M_{\mu-\frac{1}{2}}\rightarrow M_{\mu}),$$
$$W=ker(G_{-\frac{1}{2}}:M_{\mu+\frac{1}{2}}\rightarrow
M_{\mu})\cap
ker(G_{-\frac{1}{2}}G_{\frac{3}{2}}G_{-\frac{3}{2}}:M_{\mu-\frac{1}{2}}\rightarrow
M_{\mu})\cap
ker(L_{1}G_{-\frac{3}{2}}:M_{\mu-\frac{1}{2}}\rightarrow M_{\mu}).$$
Since $\mathrm{dim}M_{\mu-\frac{1}{2}}=\infty$ and
$\mathrm{dim}M_{\mu}<\infty$, $V$ is a vector subspace of finite
codimension in $M_{\mu-\frac{1}{2}}.$ Since
$\mathrm{dim}M_{\mu+\frac{1}{2}}=\infty$ and
$\mathrm{dim}M_{\mu}<\infty$, $W$ is a vector subspace of finite
codimension in $M_{\mu+\frac{1}{2}}.$ In order not to get a direct
contradiction using the Principal Fact, we assume that
$G_{\frac{3}{2}}v\neq 0$ for all $0\neq v\in V$ and
$G_{-\frac{3}{2}}w\neq 0$ for all $0\neq w\in W.$ Then
$\mathrm{dim}G_{\frac{3}{2}}V=\infty$ and
$\mathrm{dim}G_{-\frac{3}{2}}W=\infty.$

\vskip 3mm\noindent{\bf{Claim}} $\mu\neq\pm 1.$ Moreover, for any
$v\in V, w\in W,$ we have
$$G_{\frac{1}{2}}G_{-\frac{1}{2}}G_{\frac{3}{2}}v=\tau G_{\frac{3}{2}}v,\eqno(5)$$
$$G_{-\frac{1}{2}}G_{\frac{1}{2}}G_{-\frac{3}{2}}w=\tau^{'}G_{-\frac{3}{2}}w,\eqno(6)$$
where $\tau=\frac{(\mu+1)(2\mu-1)}{\mu-1}$ and
$\tau^{'}=\frac{(\mu-1)(2\mu+1)}{\mu+1}.$

\vskip 3mm\noindent{\bf{Proof\ of\ the\ Claim.}} It can be checked
directly that
$$L_{-1}G_{-\frac{1}{2}}(\frac{1}{2}L_{1}G_{\frac{1}{2}}-G_{\frac{3}{2}})
\equiv2L_0^{2}-3L_0-L_0G_{\frac{1}{2}}G_{-\frac{1}{2}}
+2G_{\frac{1}{2}}G_{-\frac{1}{2}}\ \mathrm{mod}\
U(NS)L_{-1},\eqno(7)$$ and
$$L_{1}G_{\frac{1}{2}}(\frac{1}{2}L_{-1}G_{-\frac{1}{2}}+G_{-\frac{3}{2}})
\equiv-2L_0^{2}-3L_0+L_0G_{-\frac{1}{2}}G_{\frac{1}{2}}
+2G_{-\frac{1}{2}}G_{\frac{1}{2}}\ \mathrm{mod}\
U(NS)L_{1}.\eqno(8)$$

For any $0\neq v\in V,$ by Lemma 2(i) and (7), we have
$$2L_0^{2}G_{\frac{3}{2}}v-3L_0G_{\frac{3}{2}}v-L_0G_{\frac{1}{2}}G_{-\frac{1}{2}}G_{\frac{3}{2}}v
+2G_{\frac{1}{2}}G_{-\frac{1}{2}}G_{\frac{3}{2}}v=0.$$ Then
$$(2(\mu+1)^{2}-3(\mu+1))G_{\frac{3}{2}}v-(\mu-1)G_{\frac{1}{2}}G_{-\frac{1}{2}}G_{\frac{3}{2}}v=0.$$
If $\mu=1,$ then $G_{\frac{3}{2}}v=0,$ a contradiction. Thus
$\mu\neq 1,$ and (4) holds.

For any $0\neq w\in W,$ by Lemma 2(ii) and (8), we have
$$-2L_0^{2}G_{-\frac{3}{2}}w-3L_0G_{-\frac{3}{2}}w+L_0G_{-\frac{1}{2}}G_{\frac{1}{2}}G_{-\frac{3}{2}}w
+2G_{-\frac{1}{2}}G_{\frac{1}{2}}G_{-\frac{3}{2}}w=0.$$ Then
$$(2(\mu-1)^{2}+3(\mu-1))G_{-\frac{3}{2}}w-(\mu+1)G_{-\frac{1}{2}}G_{\frac{1}{2}}G_{-\frac{3}{2}}w=0.$$
If $\mu=-1,$ then $G_{-\frac{3}{2}}w=0,$ a contradiction. Thus
$\mu\neq -1,$ and (6) holds. This completes the proof of the Claim.

\vskip3mm By (4) and (6), we see that
$G_{\frac{1}{2}}G_{-\frac{1}{2}}G_{\frac{3}{2}}v=0$ if and only if
$\mu=\frac{1}{2}$,
$G_{-\frac{1}{2}}G_{\frac{1}{2}}G_{-\frac{3}{2}}w=0$ if and only if
$\mu=-\frac{1}{2}$. There only two cases may occur: $\mu\notin
\{1,-1,\frac{1}{2}\}$ or $\mu\notin \{1,-1,-\frac{1}{2}\}.$ Because
of the symmetry of our situation, to complete the proof of the
Theorem, it is enough to show that a contradiction can be derived
when $\mu\notin \{1,-1,\frac{1}{2}\}.$

Suppose $\mu\notin \{1,-1,\frac{1}{2}\},$ then $\tau\neq0$ and
$G_{-\frac{1}{2}}G_{\frac{3}{2}}v\neq0$ for any $0\neq v\in V.$ Set
$$S=G_{-\frac{1}{2}}G_{\frac{3}{2}}V\cap W.$$
Since both
$G_{-\frac{1}{2}}G_{\frac{3}{2}}V$ and $W$ have finite codimension
in $M_{\mu+\frac{1}{2}},$ we have $\mathrm{dim}S=\infty.$

Note that
$$G_{-\frac{3}{2}}(\frac{1}{2}L_1G_{\frac{1}{2}}-G_{\frac{3}{2}})\equiv
-G_{\frac{1}{2}}G_{-\frac{1}{2}}-\frac{1}{2}L_1G_{\frac{1}{2}}G_{-\frac{3}{2}}+G_{\frac{3}{2}}G_{-\frac{3}{2}}-\frac{2}{3}c\
\mathrm{mod}\ U(NS)L_{-1}.$$ For any $v\in V,$ by Lemma 2(i), (4)
and the definition of $V$, we have
$$G_{\frac{3}{2}}G_{-\frac{3}{2}}G_{\frac{3}{2}}v=G_{\frac{1}{2}}G_{-\frac{1}{2}}G_{\frac{3}{2}}v+\frac{2}{3}cG_{\frac{3}{2}}v
=p G_{\frac{3}{2}}v, \ \textrm{for some} \ p\in
\mathbb{C}.\eqno(9)$$

Similarly, For any $w\in W,$ by
$$G_{\frac{3}{2}}(\frac{1}{2}L_{-1}G_{-\frac{1}{2}}+G_{-\frac{3}{2}})\equiv
G_{-\frac{1}{2}}G_{\frac{1}{2}}-\frac{1}{2}L_{-1}G_{-\frac{1}{2}}G_{\frac{3}{2}}-G_{-\frac{3}{2}}G_{\frac{3}{2}}+\frac{2}{3}c\
\mathrm{mod}\ U(NS)L_{1},$$ Lemma 2(ii) ,  (6) and the definition of
$W,$ we have

$$G_{-\frac{3}{2}}G_{\frac{3}{2}}G_{-\frac{3}{2}}w=G_{-\frac{1}{2}}G_{\frac{1}{2}}G_{-\frac{3}{2}}w+\frac{2}{3}cG_{-\frac{3}{2}}w
=q G_{-\frac{3}{2}}w, \ \textrm{for some} \ q\in
\mathbb{C}.\eqno(10)$$ Choose $v\in V$ such that $0\neq
G_{-\frac{1}{2}}G_{\frac{3}{2}}v\in S.$ Set $$x=G_{\frac{3}{2}}v,
y=G_{-\frac{1}{2}}x, h=G_{-\frac{3}{2}}x,
z=G_{-\frac{3}{2}}y.\eqno(11)$$  By the definitions of $W,S$ and
(10), we have
$$G_{-\frac{1}{2}}(G_{\frac{3}{2}}G_{-\frac{3}{2}}y-qy)=0-0=0$$
and
$$G_{-\frac{3}{2}}(G_{\frac{3}{2}}G_{-\frac{3}{2}}y-qy)=G_{-\frac{3}{2}}G_{\frac{3}{2}}G_{-\frac{3}{2}}y-qG_{-\frac{3}{2}}y=0.$$
So $G_{\frac{3}{2}}G_{-\frac{3}{2}}y-qy=0$ by the Principal Fact,
and we get the formula:
$$G_{\frac{3}{2}}z=qy. \eqno(12)$$
Moreover, we have
$$G_{\frac{1}{2}}z=G_{\frac{1}{2}}G_{-\frac{3}{2}}y=(2L_{-1}-G_{-\frac{3}{2}}G_{\frac{1}{2}})y
    =2G_{-\frac{1}{2}}^2y-G_{-\frac{3}{2}}G_{\frac{1}{2}}G_{-\frac{1}{2}}G_{\frac{3}{2}}v=-\tau G_{-\frac{3}{2}}x=-\tau h,$$
i.e.,
$$G_{\frac{1}{2}}z=-\tau h. \eqno(13)$$

Let $U_+$ and $U_{-}$ denote the subalgebras of $U(NS),$ generated
by $G_{\frac{1}{2}}, G_{\frac{3}{2}}$ and $G_{-\frac{1}{2}},
G_{-\frac{3}{2}},$  respectively.  We want to prove that the
following vector space
$$N=U_{+}x\oplus U_{+}y\oplus U_{-}z\oplus U_{-}h$$
is a proper weight submodule of $M,$ which will derive a
contradiction, as desired.

Since $x,y,z$ and $h$ are all eigenvectors for $L_0$, we see that
$N$ decomposes into a direct sum of weight spaces which are
obviously finite-dimensional. It remains to show that $N$ is stable
under the action of the following four operators: $G_{\frac{1}{2}},
G_{\frac{3}{2}}, G_{-\frac{1}{2}}$ and $G_{-\frac{3}{2}}$. That
$G_{\frac{i}{2}}U_{+}x\subset U_{+}x,$ $G_{\frac{i}{2}}U_{+}y\subset
U_{+}y,$ $G_{-\frac{i}{2}}U_{-}h\subset U_{-}h$ and
$G_{-\frac{i}{2}}U_{-}z\subset U_{-}z$ is clear for $i=1, 3.$

For any $a\in U_{+}, a^{'}\in U_{-},$  there exists $a_{i,j},
b_{i,j}, c_{i,j}\in U_{+}$ and $a_{i,j}^{'}, b_{i,j}^{'},
c_{i,j}^{'}\in U_{-}$ such that
$$G_{-\frac{1}{2}}a=aG_{-\frac{1}{2}}+\sum_{i,j}a_{i,j}L_0^{i}c^{j},$$
$$G_{-\frac{3}{2}}a=aG_{-\frac{3}{2}}+\sum_{i,j}a_{i,j}L_0^{i}c^{j}
+\sum_{i,j}b_{i,j}L_0^{i}c^{j}L_{-1}+\sum_{i,j}c_{i,j}L_0^{i}c^{j}G_{-\frac{1}{2}},$$
$$G_{\frac{1}{2}}a^{'}=a^{'}G_{\frac{1}{2}}+\sum_{i,j}a_{i,j}^{'}L_0^{i}c^{j},$$
$$G_{\frac{3}{2}}a^{'}=a^{'}G_{\frac{3}{2}}+\sum_{i,j}a_{i,j}^{'}L_0^{i}c^{j}
+\sum_{i,j}b_{i,j}^{'}L_0^{i}c^{j}L_{1}+\sum_{i,j}c_{i,j}^{'}L_0^{i}c^{j}G_{\frac{1}{2}}.$$
Thus, to show $G_{-\frac{i}{2}}U_{+}x,G_{-\frac{i}{2}}U_{+}y,
G_{\frac{i}{2}}U_{-}h, G_{\frac{i}{2}}U_{-}z\subset N,$ we need only
to show that $G_{-\frac{i}{2}}x, G_{-\frac{i}{2}}y,
G_{\frac{i}{2}}h,G_{\frac{i}{2}}z,G_{-\frac{i}{2}}G_{\frac{j}{2}}x,
G_{-\frac{i}{2}}G_{\frac{j}{2}}y,
G_{\frac{i}{2}}G_{-\frac{j}{2}}h,G_{\frac{i}{2}}G_{-\frac{j}{2}}z\in
N$ for $i,j=1,3.$

Now we can check the following formulas one by one via the
definitions of $V, W, S$ and (5)-(13):
$$G_{-\frac{1}{2}}x=y,$$
$$G_{-\frac{1}{2}}y=0,$$
$$G_{-\frac{3}{2}}x=h,$$
$$G_{-\frac{3}{2}}y=z,$$
$$G_{\frac{1}{2}}h=G_{\frac{1}{2}}G_{-\frac{3}{2}}G_{\frac{3}{2}}v=0,$$
$$G_{\frac{1}{2}}z=-\tau h,$$
$$G_{\frac{3}{2}}h=G_{\frac{3}{2}}G_{-\frac{3}{2}}G_{\frac{3}{2}}v=pG_{\frac{3}{2}}v=px,$$
$$G_{\frac{3}{2}}z=qy,$$
$$G_{\frac{1}{2}}G_{-\frac{3}{2}}x=0,$$
$$G_{\frac{3}{2}}G_{-\frac{3}{2}}x=px,$$
$$G_{\frac{1}{2}}G_{-\frac{3}{2}}y=-\tau h,$$
$$G_{\frac{3}{2}}G_{-\frac{3}{2}}y=qy,$$
$$G_{\frac{1}{2}}G_{\frac{1}{2}}G_{-\frac{3}{2}}y=0,$$
$$G_{\frac{3}{2}}G_{\frac{1}{2}}G_{-\frac{3}{2}}y=-\tau px,$$
$$L_{-1}x=G_{-\frac{1}{2}}^2x=0,$$
$$L_{-1}y=G_{-\frac{1}{2}}^2y=0,$$
$$G_{-\frac{1}{2}}G_{\frac{3}{2}}h=py,$$
$$G_{-\frac{3}{2}}G_{\frac{3}{2}}h=ph,$$
$$G_{-\frac{1}{2}}G_{\frac{3}{2}}z=0,$$
$$G_{-\frac{3}{2}}G_{\frac{3}{2}}z=qz,$$
$$G_{-\frac{1}{2}}G_{-\frac{1}{2}}G_{\frac{3}{2}}h=0,$$
$$G_{-\frac{3}{2}}G_{-\frac{1}{2}}G_{\frac{3}{2}}h=pz,$$
$$L_1h=G_{\frac{1}{2}}^2h=0,$$
$$L_1z=G_{\frac{1}{2}}^2z=0.$$
This completes the proof of Lemma 6 and then of Theorem
1.\hfill$\Box$

\end{document}